\documentclass[10pt,twoside]{article}
\usepackage{graphics}
\usepackage{amsmath, amsfonts, amsthm, amssymb, amscd,a4wide}
\usepackage{longtable}
\usepackage{booktabs}
\usepackage{hyperref}
\usepackage{url}
\usepackage{mathrsfs} 
\usepackage{epsfig, subfig}
\usepackage{caption}
\usepackage{multicol}
\usepackage{slashed}
\usepackage{xcolor}
\usepackage{hyperref}
\hypersetup{linktoc = all}                
\hypersetup{hidelinks}
\hypersetup{bookmarksnumbered}
\numberwithin{equation}{section}
\numberwithin{table}{section}

\newcommand \trans {\text{tran}}
\newcommand \per {\text{per}}

\newcommand \erf {\text{erf}}

\newcommand \bse {\begin{subequations}}
\newcommand \ese {\end{subequations}}

\usepackage[most]{tcolorbox}

\usepackage[most]{tcolorbox} 
\newenvironment{Figure}
  {\par\medskip\noindent\minipage{\linewidth}}
  {\endminipage\par\medskip}

\newtcolorbox{fini}{%
     enhanced, breakable, size=minimal, colframe=white, parbox=false, after={\par\vspace{2\baselineskip}}, 
     before upper={\indent}, colback=white, 
     overlay = {\draw[line width=2pt] (frame.north east) -|
                       ([xshift=3mm]frame.east)|-(frame.south east);},
     overlay first={\draw[line width=2pt] (frame.north east) -|
                           ([xshift=3mm]frame.south east);},
     overlay middle={\draw[line width=2pt] ([xshift=3mm]frame.north east) -- 
                              ([xshift=3mm]frame.south east);},
     overlay last={\draw[line width=2pt] ([xshift=3mm]frame.north east)|-
                          (frame.south east);},
}
\newtcolorbox{adone}{%
     enhanced, breakable, size=minimal, colframe=white, parbox=false, after={\par\vspace{2\baselineskip}}, 
     before upper={\indent}, colback=white, 
     overlay = {\draw[densely dotted, line width=2pt] (frame.north east) -|
                       ([xshift=3mm]frame.east)|-(frame.south east);},
     overlay first={\draw[densely dotted, line width=2pt] (frame.north east) -|
                           ([xshift=3mm]frame.south east);},
     overlay middle={\draw[densely dotted, line width=2pt] ([xshift=3mm]frame.north east) -- 
                              ([xshift=3mm]frame.south east);},
     overlay last={\draw[densely dotted, line width=2pt] ([xshift=3mm]frame.north east)|-
                          (frame.south east);},
}
\newtcolorbox{plfok}{%
     enhanced, breakable, size=minimal, colframe=white, parbox=false, after={\par\vspace{2\baselineskip}}, 
     before upper={\indent}, colback=white, 
     overlay = {\draw[dashed, line width=2pt] (frame.north east) -|
                       ([xshift=3mm]frame.east)|-(frame.south east);},
     overlay first={\draw[dashed, line width=2pt] (frame.north east) -|
                           ([xshift=3mm]frame.south east);},
     overlay middle={\draw[dashed, line width=2pt] ([xshift=3mm]frame.north east) -- 
                              ([xshift=3mm]frame.south east);},
     overlay last={\draw[dashed, line width=2pt] ([xshift=3mm]frame.north east)|-
                          (frame.south east);},
}

\newcommand \Lbf {\mathbf L} 

\newcommand \EEE E

\newcommand \bei {\begin{itemize}}
\newcommand \eei {\end{itemize}}

\newcommand \Span {\text{Span}}

\newcommand \be   {\begin{equation}}
\newcommand \bel {\be\label}
\newcommand \ee   {\end{equation}}

\newcommand \supp {{\text{supp}}}

\newcommand \ZZ    {\mathbb{Z}}
\newcommand \RR    {\mathbb R}

\newcommand \Pcal    {\mathcal{P}}

\newcommand \Hcal    {\mathcal H}

\newcommand \Lcal    {\mathcal{L}}

\newcommand \del   {\partial}

\newcommand \la         \langle
\newcommand \ra     \rangle

\newcommand \RD {{{\mathbb R}^D}}

\begin{document}

\title{The Transport-based Mesh-free Method (TMM): a review} 

\author{Philippe G. LeFloch\footnote{Laboratoire Jacques-Louis Lions, Centre National de la Recherche Scientifique, Sorbonne Universit\'e, 4 Place Jussieu, 75252 Paris, France. Email: {\tt contact@philippelefloch.org}.}
\, 
 and Jean-Marc Mercier\footnote{MPG-Partners, 136 Boulevard Haussmann, 75008 Paris, France. Email: {\tt jean-marc.mercier@mpg-partners.com.}}
}

\date{October 2019}

\maketitle

\begin{abstract}   
We review a numerical technique, referred to as the Transport-based Mesh-free Method (TMM), and we discuss its applications. We recently introduced this method from a numerical standpoint and investigated the accuracy of integration formulas based on the Monte-Carlo methodology: quantitative error bounds were discussed and, in this short note, we outline the main ideas of our approach. 
The techniques of transportation and reproducing kernels lead us to a very efficient methodology for numerical simulations in many practical applications, and provide some light on the methods used by the artificial intelligence community. For applications in the finance industry, our method allows us to compute many types of risk measures with accurate and fast algorithms. We propose theoretical arguments as well as extensive numerical tests in order to justify sharp convergence rates, leading to rather optimal computational times. 
Cases arising in finance applications support our claims and, finally, the problem of the curse of dimensionality in finance is briefly discussed. 
\end{abstract}

\begin{multicols}{2}
 
\section{Introduction}
\label{section:introd}

Relying on our recent papers \cite{PLF-JMM-2}--\cite{PLF-JMM-4}, we present and discuss here a numerical technique, that we refer to as the Transport-based Mesh-free Method (TMM), which is of direct interest in numerical simulations. 
Our method is mesh-free (cf.~for instance \cite{FGE2,Wendland-book}) and somewhat similar to a Lagrangian mesh-free method. Importantly, our method can handle transport as well as diffusive terms and was introduced first in \cite{PLF-JMM-2}. 

Our motivation was to reduce as much as possible the algorithmic burden of solving partial differential equations (PDEs) especially for problems in large dimensions, met for instance in mathematical finance and machine learning. Computational times reflect, in a concrete manner, the complexity of an algorithm. For PDEs solvers, the algorithmic complexity can be measured by establishing suitable error estimates. For our TMM approach, in \cite{PLF-JMM-3}  we were able to establish some Monte-Carlo type error estimates, at least via heuristic arguments, as we outline below in Section~\ref{SMCE}. 

This allowed to perform a precise error analysis of this method in \cite{PLF-JMM-4}, which is outlined in Section~\ref{TMMF}.  Finally, in Section~\ref{COD} we discuss the limitations coming from the curse of dimensionality for applications to finance.

The TMM methodology has wide applications in mathematical finance, since it allows one to compute almost any risk measures, quite accurately and with a fast algorithm. A risk measure is here understood as a price, future prices, future sensitivities, Value at Risk (VaR), or Counterparty Value Adjustment (CVA), and may concern a simple asset, a complex derivative or an investment strategy, 
as well as a big portfolio of such instruments; they can be 
written on any number of underlyings, themselves depending on any Markov-type stochastic processes.  

The proposed method was extensively tested in mathematical finance ones; 
see \cite{PLF-JMM-2}-\cite{PLF-JMM-4} as well as  \cite{JMM-SM} for a business case in asset and liability management using the so-called Libor market model \cite{BGM1997}. Another business case, for front-office equity derivatives, was treated using this method: it consists in computing metrics for specific customer needs for big portfolios of autocalls, that are useful for pre-sales purposes. The modeling of shares uses the Buelher dividend models \cite{Bu}, and the algorithm described in \cite{PLF-JMM-2} for local volatility calibration. For an application to nonlinear propagation, see \cite{PLF-JMM-1}-\cite{PLF-JMM-4}.


\section{Monte-Carlo-type strategy}\label{SMCE}

We postpone the discussion of earlier references at end of this section and outline now our strategy for deriving a priori error estimates on multidimensional integrals. 
One of our task is to investigate the validity of Monte-Carlo-type error estimates of the form
\bel{MCE}
   \Big| \int_{\RR^D} \varphi(x)d\mu - \frac{1}{N} \hskip-.2cm
\sum_{1 \leq n \leq N} \varphi(y^n) \Big| \le E_K(Y,N,D) \|\varphi\|_{\Hcal_K}.
\ee
Here, \(\mu \in \Pcal(\RR^D)\) is a probability measure (whose support $\supp(\mu)$ must be convex) and $Y = (y^1,\ldots,y^N)$ is a set of $N$ distinct points in 
$\RR^{N}$. 
We have denoted here by $\Hcal_K$ a kernel-based Hilbert space depending upon the choice of an admissible kernel $K$, that is, a continuous and symmetric function $K: (x,y) \in \RD \times \RD \mapsto \RR$ with $K(x,y) = K(y,x)$. Admissibility means that the matrix 
\bel{ADK}
   K(Y,Y):= \big( K(y^n,y^m) \big)_{1 \leq n,m \leq N}
\ee
is symmetric positive-definite for any choice of $Y$. The function space $\Hcal_K$ is sometimes called a reproductible Hilbert kernel space (RHKS) or a native space. The terminology is a little bit confusing, since \eqref{MCE} is relevant for Hilbert spaces as as their generalization to the corresponding Banach spaces. 
The Hilbert space of interest here consists of all linear combinations of the functions $K(x, \cdot)$ (parametrized by $x \in \RR^D$), that is, 
\bel{HK}
\Hcal_K
:=  \Span \big\{K(\cdot, x) \, / \, x \in \RR^D \big\}, 
\ee
endowed with a norm induced by a scalar product defined such that
$$
\big\la K(\cdot, x), K(\cdot, y) \big\ra_{\Hcal_K} 
= K(x,y), \qquad x, y \in \RR^D. 
$$

In \eqref{MCE}, the function $E_K(Y, N, D)$ is referred to as the {\sl discrepancy error function} and can be expressed as:
\bel{EE}
\begin{aligned}
   & E_K(Y,N,D)^2
 = \iint_{\RR^{D} \times \RR^{D}} K(x,y)d\mu_xd\mu_y 
\\
& + \frac{1}{N^2} \hskip-.2cm
\sum_{n,m=1}^N \hskip-.2cm
K(y^n,y^m)- \frac{2}{N} \sum_{1 \leq n \leq N} \int_{\RR^D} K(x,y^n)d\mu_x.
\end{aligned}
\ee
Observe that this error function can be readily approximated by using, for instance, a direct Monte-Carlo approach. We assuming that $K$ is integrable for the measure $\mu$ with respect to both variables. Then we say that a sequence $\overline Y$ 
is \textbf{sharp discrepancy sequence} if it achieves the global minimum of the functional, that is, 
\bel{SDS}
  \overline{Y} = \arg \inf_{Y \in \RR^{N \times D}} E_K(Y,N,D),
\ee
and we denote the minimum by 
\bel{SDE}
E_K(N,D) = E_K(\overline{Y},N,D). 
\ee
Of course, it practice, we need achieve exactly the minimum and $E_K(\overline{Y},N,D)$ serves as our {\bf error discrepancy bound} when $\overline Y$ is our numerical solution. 

The overall construction is as follows.  
To any admissible kernel $K$ we associate the function space $\Hcal_K$ in \eqref{HK}, within which the accuracy of a numerical approximation formula 
can be evaluated by computing the error function \eqref{EE}. 
In order to optimize the convergence rate arising in \eqref{MCE}, we should choose the points in order to achieve \eqref{SDS}. Moreover, whenever we are able to compute (and at least estimate) the discrepancy error \eqref{SDE}, then we have a method for evaluating {\sl quantitatively} the accuracy of our approximation.

We now list several important classes of admissible kernels:
\begin{itemize}

\item {\bf Translation-invariant kernels,} by definition, have the form $K(x,y) = \nu(x-y)$, where $\nu$ is a function whose Fourier transform is a probability measure, that is,
 $\widehat{\nu} \in \Pcal(\RR^D)$ (thanks to Bochner theorem). Among them, one can consider the important class of {\bf radially-symmetric kernels} $\nu(|x-y|)$, including kernels generating the standard Sobolev spaces. 
We emphasize that such kernels are \textit{not localized} in the sense that $\nu(x-y)$ fails to be in (possibly weighted) $L^1(\RR^D \times \RR^D)$.

\item {\bf Zonal  kernel} \cite{MA} or {\bf power series kernels} \cite{ZB}, by definition, have the form $K(x,y) = F(<x,y>)$, where $<\cdot, \cdot>$ denotes the Euclidian scalar product in $\RR^D$ and $F: \RR \to \RR$ is fixed and is called an activation function. Such kernels are used by the artificial intelligence community, together with {\bf convolutional kernels,} which are translation-invariant kernels of the form $K(x,y) = (\phi \ast \phi)(x-y)$ (where $\ast$ denotes the convolution operator)s.
\end{itemize}


Throughout, we are given a convex and open set $\Omega \subset \RD$ which is assumed to have a piecewise smooth boundary and, typically, we will take $[0,1]^D$. 
We observe that, using a transportation argument, it is sufficient take in \eqref{MCE}  the Lebesgue measure $\mu = dx_\Omega$ on $\Omega$. 
Namely, 
if $S:\Omega \mapsto \RR^D$ is a transport map for a general measure $\mu$, that is, 
the unique map satisfying $\int_{\RR^D} \varphi d\mu = \int_{\Omega} (\varphi \circ S)dx$ for any continuous $\varphi \in L_\mu^1(\RD)$, together with $S = \nabla h$, $h$ convex and $\nabla$ the gradient operator. Indeed, using such a map, \eqref{MCE} can be written as
\bel{MCEL}
\aligned
	&   \Big| \int_{\Omega} (\varphi \circ S)dx - \frac{1}{N}\sum_{1 \leq n \leq N} (\varphi \circ S)(x^n) \Big|
\\
& \le E_K\big(X,N,D) \|\varphi\|_{\Hcal_K},
\endaligned
\ee
with $y^n = S(x^n)$ and $\varphi \circ S$ denoting the composition of two functions. 

Let us briefly review some of the earlier literature about the estimate \eqref{MCE} arising in approximation theory. One of the most used integration method is the direct Monte-Carlo method, and is based on i.i.d. sequences $Y$; 
with suitable statistical arguments and applying the law of large numbers, one can estimate $E_K(Y,N,D) \sim \frac{1}{\sqrt{N}}$ with a variance-type norm, that is, the space $\Hcal_K$ is replaced by $L_\mu^2(\RD)$. Low-discrepancy sequences (see \cite{Ni} and the references therein) and Sobol sequences lead to estimates
in the bounded variation space $BV([0,1]^D)$ and, specifically, it is expected that $E_K(N,D) \lesssim \frac{\ln(N)^{D-1}}{N}$
---a bound referred to as the \textit{Koksma--Hlawka conjecture}. 

Many other estimates of this type are available in the 
literature concerned with 
wavelets, quantification, neural networks. 
Notably, for mesh-free methods, Wendland and followers derived error estimates with radial basis functions
in the 90's; see for instance \cite{NWW}. Our contribution (see next section)
is a systematic study of the discrepancy error function for a variety of  admissible kernels.

\section{Kernel-based estimates}

A kernel is usually chosen and adapted to a specific application. Once chosen, the accuracy of the method will eventually depend, as described in the previous section, upon our ability to solve the minimization problem \eqref{SDS}-\eqref{SDE}. 

For applications in finance, we carefully designed kernels adapted to several important requirements. For instance, we present here kernels based on the tensor-based Matern kernel, which is adapted to spaces of functions $\varphi =\sum_{0<n_1< \ldots< n_k \le D } \varphi_{n_1,\ldots,n_k}(x_{n_1},\ldots,x_{n_k})$. In particular, this choice appears to be well-adapted to describe a {\sl portfolio structure.}
 Moreover, we found it as well important to have {\sl localized kernels,} that are kernels defined on a simple set $\Omega$, such as $[0,1]^D$. We investigated 
two localization techniques:
\begin{itemize}

\item \textbf{Periodic kernels} based on a discrete lattice are motivated by the work by Cohn and Elkies \cite{CohnElkies} who studied the problem of sphere packing. Consider a family of $D$ vectors $l_1, l_2, \ldots, l_D \in \RD$ being given, and define the lattice $\Lbf:= \Big \{\sum_{1 \leq d \leq D} \alpha_d  l_d, \ \alpha_d \in \ZZ \Big \}$, and its dual lattice $\Lbf^*:= \Big\{\alpha^* \in \RD \, \big/ \,  <\alpha, \alpha^*> \in \ZZ \, \text{ for all } \alpha \in \Lbf \Big \}$. Consider any discrete function satisfying $\rho(\alpha^*) \in \ell^1(L^*)$ with $\rho(\alpha^*) \ge 0$ and $\rho(0)=1$. Then, a {\bf lattice-based kernel} is the $L$-periodic, translation-invariant kernel 
\bel{LBK}
K^\per(x,y) = \frac{1}{|C|}\sum_{\alpha^* \in \Lbf^*} \rho(\alpha^*)e^{2 i \pi <x-y,\alpha^*>},
\ee
where $|C|$ is the volume of the elementary cell $C$ defining the lattice.

\item \textbf{Transported kernels} are defined from prescribing an admissible kernel $K$, (for instance a lattice-based one),
 and a transport map $S:\Omega \mapsto \RR^D$, with $S = \nabla h$ and $h$ convex. Based on these data, we then introduce the kernel 
\bel{TK}
   K^\trans(x,y) = K(S(x),S(y)). 
\ee

\end{itemize}
In Figure \ref{Kernels} we illustrate these two localization techniques. 
We plot the tensorial version of the Mat\'ern kernel (also called exponential kernel), that is, the translation-invariant kernel $K(x,y):= \chi(x-y) = \exp(-|x-y|_1)$, in which 
$|x|_1 = \sum_{d=1}^D |x_d|$. We use here $\rho(\alpha) = \widehat\chi(\alpha)$ in \eqref{LBK} to define the lattice-based kernel, while
 the transport map in \eqref{TK} is chosen to be $S(x) = \erf(x)$. Here, $\erf$ denotes  the standard error function (i.e.~the integral of the normal distribution). 
We plots $K^\per(x,0)$ on the left-hand side
 and $K^\trans(x,0)$ on the right-hand side 
with $\Omega = [0,1]^D$ and $D=2$.

\begin{Figure}
\includegraphics[width=0.49\linewidth]{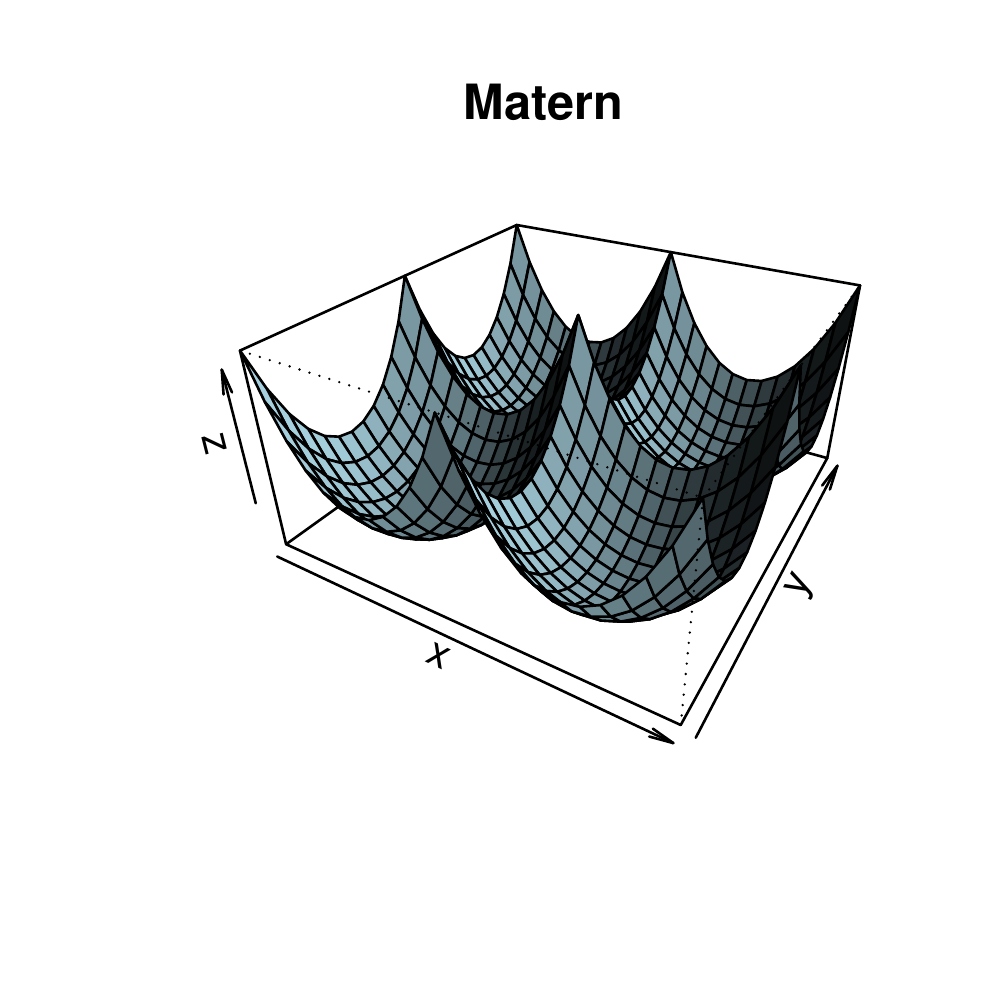} 
\includegraphics[width=0.49\linewidth]{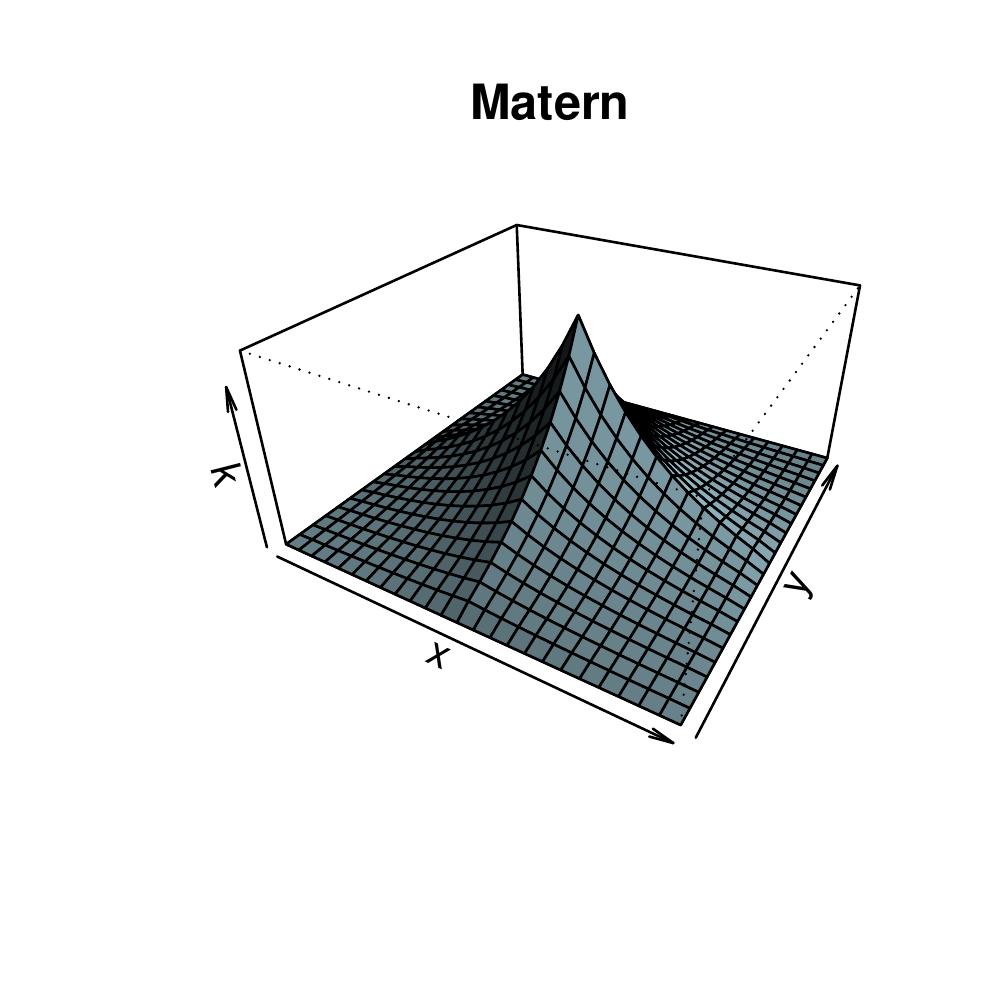} 
\captionof{figure}{Periodic/transported Mat\'ern kernels}
\label{Kernels}
\end{Figure}


For both localization techniques, we can approximate the sharp discrepancy sequences in view of \eqref{SDS}. For instance, Figure \eqref{fig:DISTRIBLATTICE} shows three distributions in the two-dimensional case: the first one is a random Mersenne Twister sequece (MT19997); 
the second one is a sequence approximating the sharp discrepancy one for the lattice-based Mat\'ern kernel (from the left-hand side of Figure \eqref{Kernels});
the third one corresponds to the transported Gaussian kernel and
is designed from the Gaussian kernel $K(x,y) = \exp(-|x-y|^2)$, that is a translation-invariant, radially-symmetric kernel, to which we applied the transport map $\erf$. 

Observe that the distribution corresponding to the Gaussian kernel 
ressembles an optimal sphere packing. 
On the other hand, the distribution associated with the Mat\'ern kernel can also be 
interpreted as an optimal packing (but not a sphere packing).

\begin{Figure}
\includegraphics[width=0.32\linewidth]{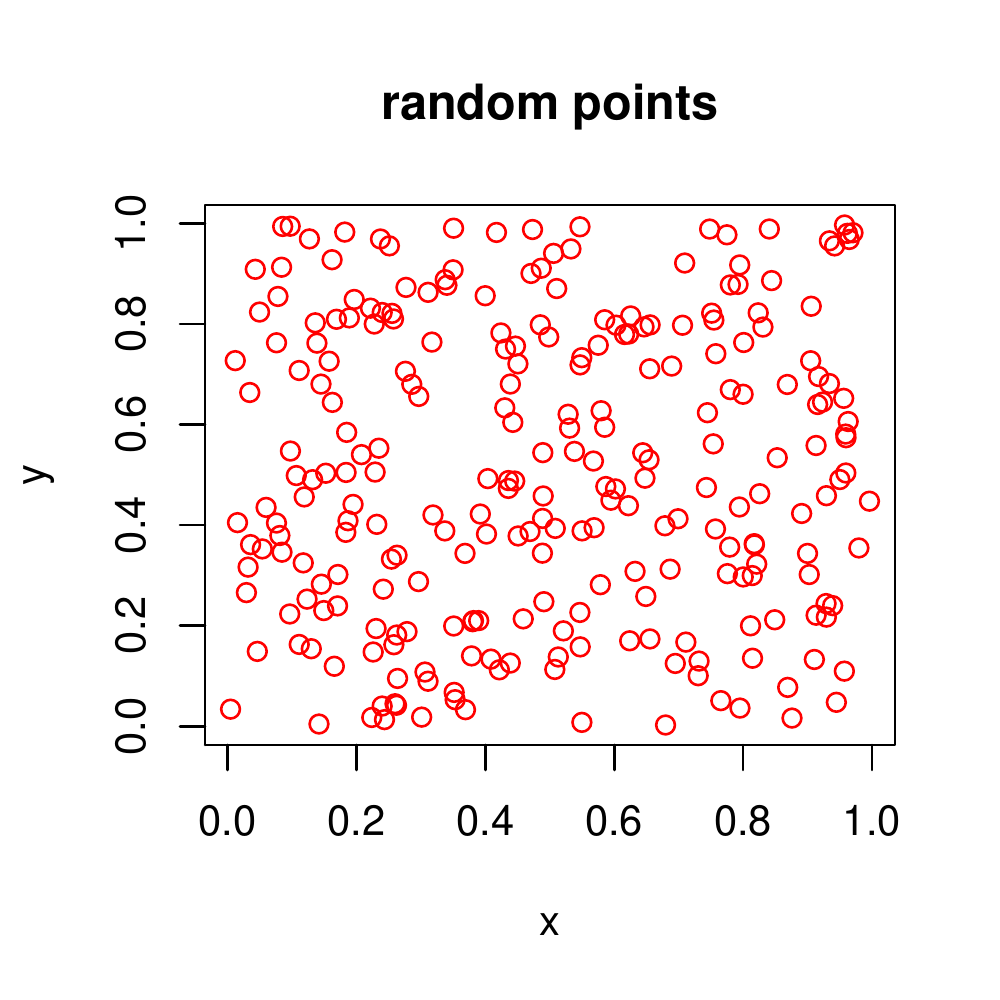} \includegraphics[width=0.32\linewidth]{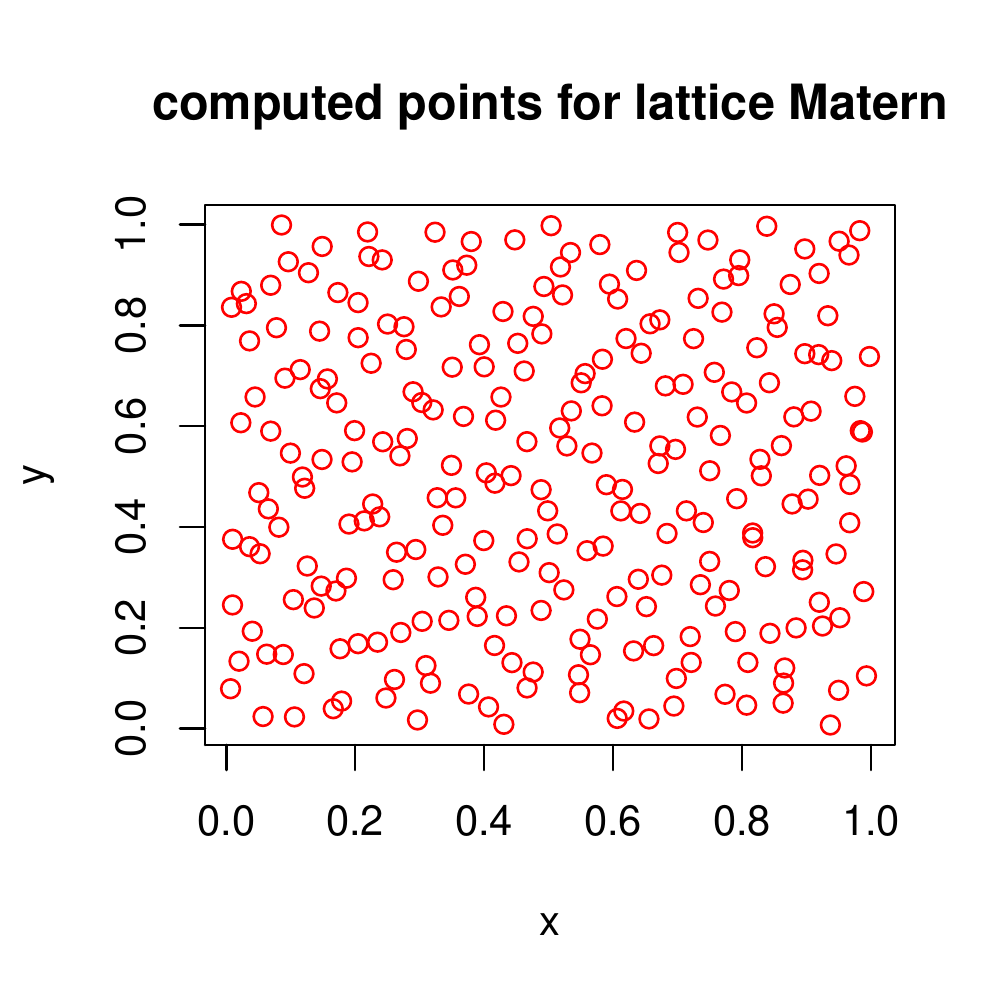} 
\includegraphics[width=0.32\linewidth]{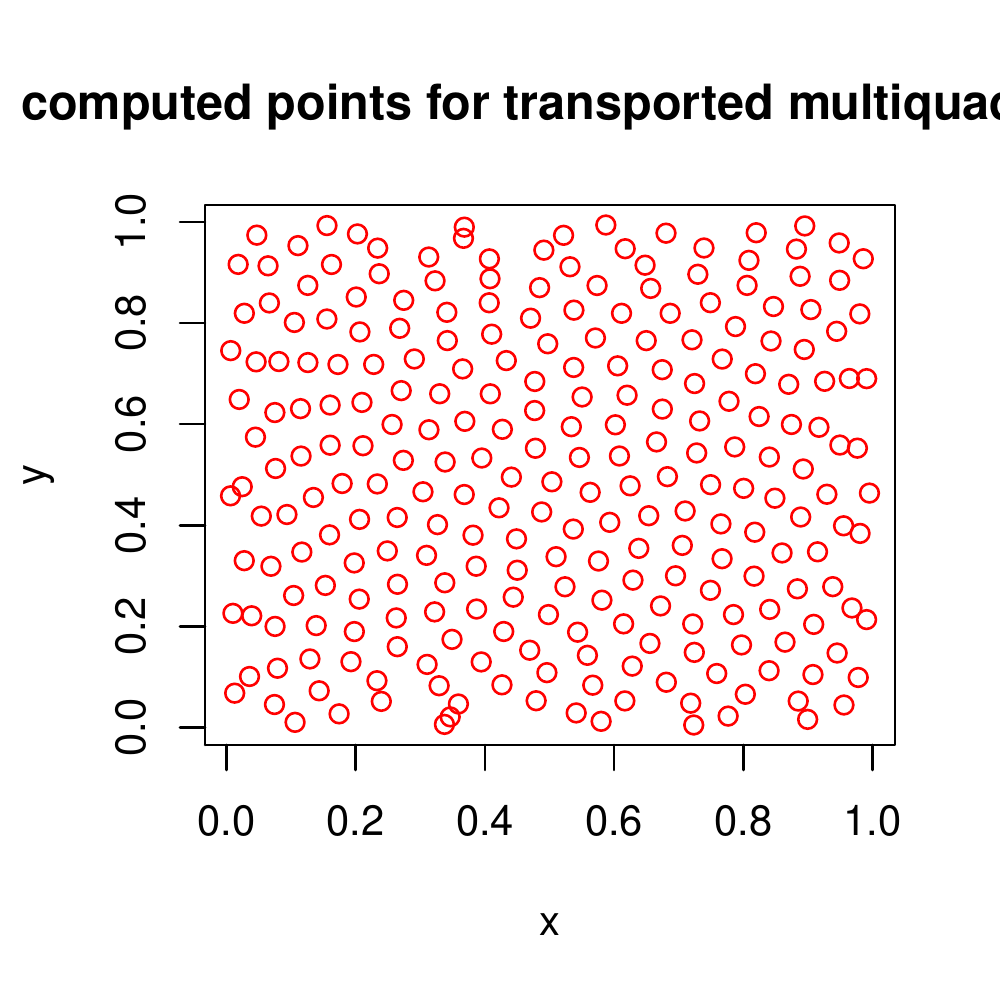}\captionof{figure}{Random/Mat\'ern/Gaussian distributions with $N=256$}\label{fig:DISTRIBLATTICE}
\end{Figure}


However, while considering lattice-based kernels \eqref{LBK}, our theoretical 
analysis can be supported by the following formula, which  provides quantitative  information on the error \eqref{SDE}: 
\bel{SEA}
    \overline{E}_K(N,D) \lesssim \sqrt{\frac{1}{N} \sum_{n > N} \phi(\alpha^{*n})},
\ee
where the ordering chosen for the lattice points $\alpha^{*n}$ is such that
the sequence $\phi(\alpha^{*n})$ is decreasing .
 This formula can be evaluated numerically or theoretically using a level-set argument; 
see \cite{PLF-JMM-3} for the details.

For instance, the following table was obtained in \cite{PLF-JMM-3} for the lattice-based Mat\'ern kernel, in which we compare \eqref{SEA} with the discrepancy error obtained from the minimization problem \eqref{SDS} for a broad range of values $N$ and dimensions $D$. As can be seen from this table, the error approximation formula \eqref{SEA} is not perfectly exact, but does 
give a good idea of the accuracy of the computed sequence. Moreover, the estimate \eqref{SEA} can be roughly approximated as $E_K(N,D) \lesssim \frac{\ln(N)^{D-1}}{N}$, hence similar to the one in the Koksma--Hlawka conjecture. 

\begin{center}
\scalebox{0.6} {
\begin{tabular}{|l||c||c||c|}
  \hline
    & D=1 & D=16 & D=128 \\
  \hline
   N=16 & 0.062 & 0.211 & 0.223 \\
   N=128 & 0.008 & 0.069 & 0.077 \\
  N=512 & 0.002 & 0.034 & 0.049 \\
  \hline
\end{tabular} \ \
\begin{tabular}{|l||c||c||c|}
  \hline
    & D=1 & D=16 & D=128 \\
  \hline
   N=16 & 0.062 & 0.288 & 0.323 \\
   N=128 & 0.008 & 0.077 & 0.105 \\
  N=512 & 0.002 & 0.034 & 0.043 \\
  \hline
\end{tabular}
}
\end{center}

\section{Main equations for finance}

Before we can outline our approach (in the next section), let us we briefly describe the equations that one solves in finance, that is, 
the Fokker-Planck and Kolmogorov equations. 
We begin with the definition of a stochastic differential equation (SDE) describing the dynamics of a Markov-type stochastic process, denoted by $t \mapsto X_t \in \RD$, i.e. 
\bel{SED}
dX_t = r(t,X_t)dt+\sigma(t,X_t)dB_t. 
\ee
Here, $B_t \in \RR^D$ denotes a $D$-dimensional, independent Brownian motion, 
while 
$r \in \RR^D$ is a prescribed vector field and $\sigma \in \RR^{D\times D}$ is a prescribed matrix-valued field. 

Denote by $\mu = \mu(t,s,x,y)$ (defined for $t\ge s$) the {\bf density probability measure} associated with $X_t$, \textit{knowing} the value $X_s = y$ at the time $s$. We recall that $\mu$ obeys the {\bf Fokker-Planck equation,} which is 
the following nonlinear partial differential equation (defined for $t \geq s$): 
\bel{FP}
	\del_t \mu - \Lcal \mu = 0, \quad \mu(t,s,x, y)|_{t=s, x=y} = \delta_{y},
\ee 
$\delta_{y}$ being the Dirac mass weighting $y$, which is a convection-diffusion equation. Moreover, the initial data is 
the Dirac mass $\delta_y$ at some point $y$, while
 the partial differential operator is  
\bel{Lcal} 
\Lcal \mu:= \nabla \cdot (r \mu) + \nabla^2 \cdot (A \mu), \quad A:= {1 \over 2} \sigma \sigma^T. 
\ee
Here, $\nabla$ denotes the gradient operator, $\nabla \cdot$ the divergence operator, and $\nabla^2:= (\del_i \del_i)_{1 \leq i,j \leq D}$ is the Hessian operator. 
We are writing here $A \cdot B$ for the scalar product associated with the Frobenius norm of matrices. We emphasize that weak solutions to \eqref{FP} defined in the sense of distributions must be considered, since the initial data is a Dirac mass.

The (vector-valued) dual of the Fokker-Planck equation is the {\bf Kolmogorov equation}, also
 known in mathematical finance as the {\bf Black and Scholes equations}.
For an unknown $\overline{P} = \overline{P}(t,x)$ with $t \leq s$ 
 reads  
\bel{KE}
\aligned
&	\del_t \overline{P} - \Lcal^* \overline{P} =0, 
\\
&
\Lcal^* \overline{P}:= - r \cdot \nabla \overline{P} + A \cdot \nabla^2 \overline{P}, 
\endaligned
\ee
and the vector-valued function $\overline{P} \in \RR^M$ models a portfolio of $M$ instruments where $M$ is typically a {\sl large} integer.
The
Kolmogorov equations \eqref{KE} are the equations of interest that one solves
in the applications to finance: they determine the so-called \textit{fair values}. 
Namely, thanks to the Feynmann-Kac theorem, a solution to the Kolmogorov equation \eqref{KE} can be interpreted to be a time-average of an expectation function, as follows: 
\bel{Ptsr}
\overline{P}(t,s,y) 
= \int_{s \le u } \mathbb{E}^{X_u}\Big[ P(t,u,X_u) | \ X_s=y\Big]
 du, 
\ee
in which $P(t,s,X)ds$ is called the \textit{payoff} of any instruments whose underlying is described by the random variable $X$.
Here, we distinguish the payoff \(P\) from its fair value, using the overline notation \(\overline{P}\). For instance,
provided \(s > t\) ($t$ being `today'), then $\mathbb{E}^{X_s} \Big[ \overline{P}(t,s,\cdot) | \ X_t=y\Big]$
is called the \textbf{forward value} of the instrument at the time $s$.

Solving the Kolmogorov equations for a given instrument allows one to compute not only its price ---which is $\overline{P}(0,y) = \overline{P}(0,y)|_{(t,x) = (0,y)}$ in the above setting---
 but also all of the \textit{fair value surface} $(t,x) \mapsto \overline{P}(t,x)$
(for all $t \ge 0$ and $x \in \RD$). This latter observation 
is important in an operational context, since all standard risk measures can be determined from the knowledge of this surface, such as 
risk measures of internal or regulatory nature, 
or optimal investment strategies: 
for instance, American exercising, or sophisticated hedging strategies based on sensitivities \cite{JMM-SM}. 


\section{TMM in finance}\label{TMMF}

Our numerical strategy, which we refer to as the transport-based mesh-less method, 
allows to solve the above two equations, namely the Fokker-Planck and the Kolmogorov equations. Here, we only outline the arguments and explain how 
quantitative error estimates are be ensured; we refer the reader to \cite{PLF-JMM-2} and \cite{PLF-JMM-4} for further details. We emphasize however that the proposed framework can be used for more general problems of hyperbolic-parabolic type, such as the Hamilton-Jacobi equations \cite{PLF-JMM-1}, Euler equations, and Navier-Stokes equations. 

\textbf{Step 1: the forward computation.} 
Consider the Fokker-Planck equation \eqref{FP} together with the Monte-Carlo-type error estimate \eqref{MCE}. Once a kernel $K$ is selected, we can apply the numerical scheme presented earlier in \cite{PLF-JMM-2}, which is 
a stable and consistent approximation of the Fokker-Planck equation \eqref{FP}
and provides an approximation of the solution $\mu$. 
This approximation is a discrete probability measure of the form 
$\frac{1}{N} \big(\delta_{y^1(t)} +\ldots + \delta_{y^N(t)}\big)$. 
 Interestingly enough, $Y(t) = (y^1, \ldots,y^N)(t)$ converges toward a sharp discrepancy sequence, in the sense defined in \eqref{SDS}. To check the accuracy of this numerical step, at each discrete time we can compute the error discrepancy \eqref{EE}. That is, we have the following error estimate for any moment of the measure $\mu$ at any time $t$: 
\bel{FPE}
\begin{aligned}
& \Big| \int_{\RR^D} \varphi(x)d\mu(t) - \frac{1}{N}\sum_{1 \leq n \leq N} \varphi(y^n(t)) \Big| \\
& \le E_K\big(Y(t),N,D\big) \|\varphi\|_{\Hcal_K}
\end{aligned}
\ee
for any continuous test-function $\varphi \in L_{\mu(t)}^1(\RD)$. Then, in view of the sharp discrepancy error \eqref{SDE}, we can compare $E_K\big(Y(t),N,D\big)$ with the exact minimum value $E_K(N,D)$ and, therefore, {\sl explicitly check the accuracy} of the numerical solution.

Let us illustrate this procedure with the (shifted) SABR model (see \cite{AKS2015} and the references therein) for a time evolution
with initial conditions \(F_0\) and \(\alpha_0\), described by the
following coupled system of stochastic differential equations:
\begin{equation}
d \Big( \begin{array}{c}
 F_t \\ 
 \alpha_t
\end{array} \Big) = \rho \Big( \begin{array}{cc}
\alpha_t (F_t+s)^\beta & 0 \\ 
0 & \nu \alpha_t
\end{array} \Big)
\Big( \begin{array}{c}
 dW_t^1 \\ 
 dW_t^2
\end{array} \Big). 
\ee
Here \(0\leq \beta \leq 1\) is a parameter representing the constant elasticity of variance
(CEV), \(\nu \geq 0\) is a constant volatility
parameter, \(W^1_t, W^2_t\) are two independent Brownian motions, and 
$\rho$  is a real-valued correlation matrix.

Consider the transported Mat\'ern kernel in Figure \ref{Kernels} (right-hand figure). Then Figure \ref{SABRden} is a plot of our approximation of the sharp discrepancy sequence for the SABR model; 
it uses \(N=200\) points with the parameters $F_0=3\%$, $\alpha_0 = 10\%$, $\nu = 10\%$, $\beta = 1$, and $\rho_{12} = \rho_{21} = 0.5$. This figure shows
 the set \(\big( y^1(t),\ldots, y^N(t)\big)\), where the \(y\)-axis represents the volatility process \(\alpha_t\) and \(x\)-axis the values of the interest rates \(F_t\).

\begin{Figure}
\includegraphics[width=0.32\linewidth]{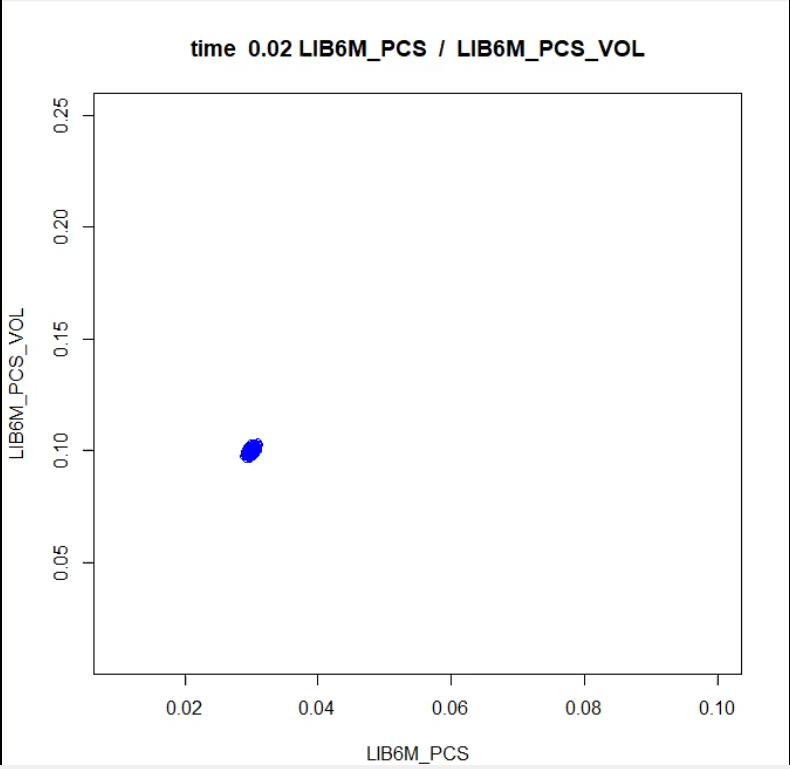} \includegraphics[width=0.32\linewidth]{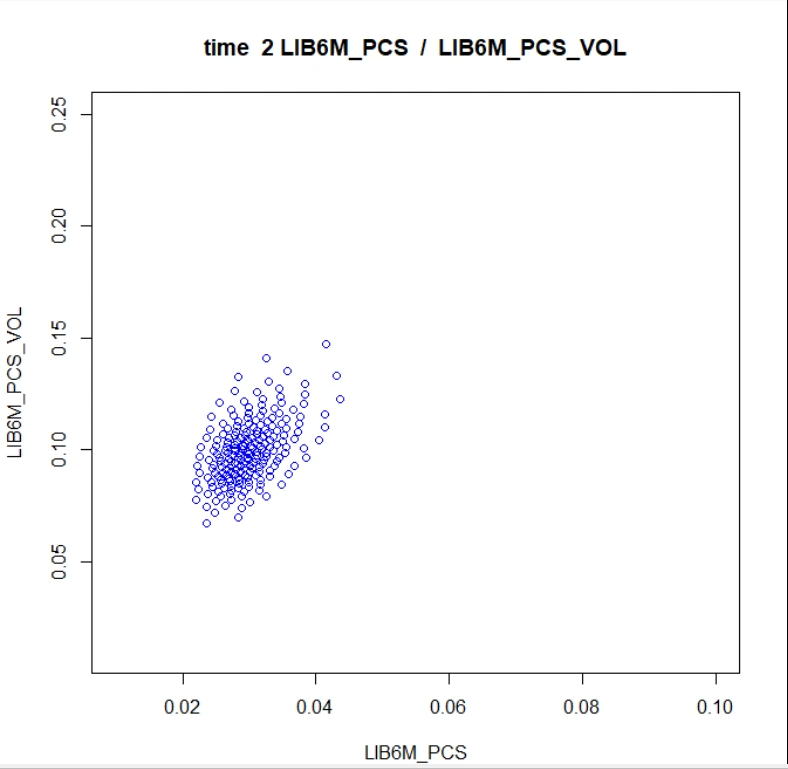} 
\includegraphics[width=0.32\linewidth]{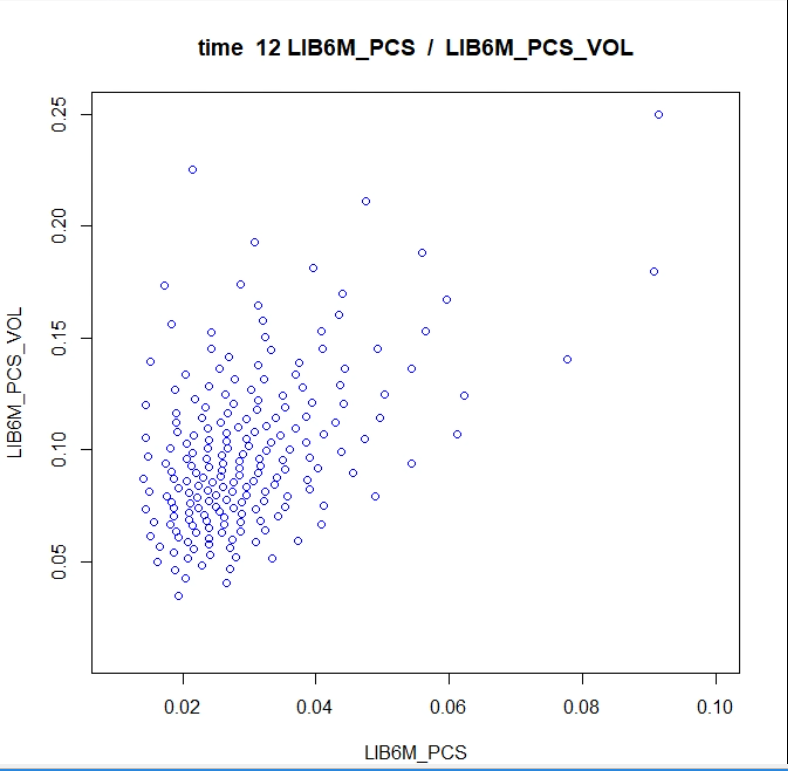}
\captionof{figure}{SABR at time 0.02, 2 and 12. N=200.}\label{SABRden}
\end{Figure}

\textbf{Step 2: the backward computation.} 
Once the sharp discrepancy sequence is computed, we are in a position to 
solve the Kolmogorov equation \eqref{KE}, using $t \mapsto y^n(t)$ (with $n= 1, \ldots, N$) as a moving transported grid, again using the numerical scheme in \cite{PLF-JMM-2}. This scheme provides us with an approximation which is consistent with the Kolmogorov equation and, in view of \eqref{FPE} we see that this scheme enjoys the error estimate 
\be
\aligned
& \Big| \int_{\RR^D} \overline{P}(t,\cdot)d\mu(t,\cdot) - \frac{1}{N}\sum_{1 \leq n \leq N} \overline{P}(t,y^n(t)) \Big| 
\\
& \le 
E_K\big(Y(t),N,D\big) \|\overline{P}(t,\cdot)\|_{\Hcal_K}.
\endaligned
\ee
Moreover, the discrete solution $t \mapsto P(t) \in \RR^{N \times M}$ approaching $\big( P(t,y^n(t)) \big)_{1 \leq n \leq N}$ is computed accordingly to
\bel{KED} 
\overline{P}(s) =  \Pi^{ (t,s)}\overline{P}(t),  
\quad \Pi^{ (t,s)}:= \Big( \pi^{ (t,s)}_{n,m} \Big)_{1 \leq n,m \leq N},
\ee
where the matrix $\Pi^{ (t,s)} \in \RR^{N \times N}$  is computed explicitly
and is nothing but the generator of the discrete counterpart of the Kolmogorov equation. 
 The matrix $\Pi^{ (t,s)}$ is interpreted in a Markov-chaining process setting as follows: 
 $\pi^{ (t,s)}_{n,m}$ is the probability that the stochastic process jumps from the sharp discrepancy state $y^n(t)$ to the sharp discrepancy state $y^m(s)$. Indeed, our numerical scheme, by construction, yields this matrix as a stochastic matrix --or a bi-stochastic matrix (i.e.~having each row and column summing to $1$) if the underlying is a martingale process. 

We point out that we can also treat a boad set of partial derivative operators and, for instance, forward sensitivities:
\bel{KEDG}
	\nabla \overline{P}(s) \text{ approximates } \nabla_y \overline{P}(t,y^n(t))_{1 \leq n \leq N}.
\ee
This allows us to compute, for instance, hedging strategies \cite{JMM-SM}. We can also treat more complex operators such as the Hessian operator or the Helmholtz-Hodge decomposition, which are important in, for instance, fluid dynamics.


\section{Remarks on the curse of dimensionality in finance}
\label{COD}

Let us emphasize that our method shed some new light on the problem of the curse of dimensionality for applications to finance. This classical problem is stated as follows: consider a stochastic process modeling several underlyings $t \mapsto X_t \in \RD$ (with $D>>1$), 
and consider a payoff of a complex option $P(t,X_t)$. In order to manage such an instrument, we would like to have some definite confidence on the numerical algorithm that we use for computing its fair values or its sensitivities.

Consider first a lattice-based periodic kernel $K$ and the formula \eqref{LBK} for which we can specify directly its Fourier coefficients $\rho(\alpha)$. In particular, using for instance the estimate \eqref{SEA}, our algorithm \eqref{KED} for the Kolmogorov equation provide an approximation at any order of accuracy $a \ge 1/2$: 
\be
\aligned
& \Big| \int_{\RR^D} \hskip-.2cm 
\overline{P}(t,\cdot)d\mu(t,\cdot) & - \frac{1}{N}\sum_{1 \leq n \leq N} \overline{P}(t,y^n(t)) \Big| 
 \\
& \le \frac{\|\overline{P}(t,\cdot)\|_{\Hcal_K}}{N^{a}}. 
\endaligned
\ee
The limit case $a = \infty$ is quite intriguing, and we can then also choose the function $\rho(\alpha) = 1$ if $\alpha = 0$, while 
$\rho(\alpha) = 0$ otherwise. With this limiting choice, the function space $\Hcal_K$ contains constant functions only and, of course, most of the `information'' on the function is lost.  However, the main point is that 
we can achieve any order of convergence at the expense of increasing the decay of 
the Fourier coefficients that determine a lattice-based kernel. This in turn defines function spaces $\Hcal_K$ of functions that are more regular as the dimension increases.

This above effect, in principle, could be problematic while managing a financial instrument that has a rather low regularity. For instance, American-type options require kernels and  modeling functions whose second-order derivatives are only measures. This is even worse for autocalls, that are functions whose first-order derivatives are signed measures. Hence, for such instruments it is very desirable to {\sl carefully quantify the numerical error} made in computing prices and derivatives. The error estimates presented in this paper can  be very helpful for this purpose.


\section{Conclusions}

In this note based on \cite{PLF-JMM-2}--\cite{PLF-JMM-4},
we presented a new analysis of  Monte-Carlo-type integration formula, which is relevant in a variety of applications and leads to sharp error estimates of practical interest.

We also presented a new numerical method, which we refer to as the Transport-based mesh-free Method, and is designed for the numerical simulations of PDEs and should be useful for a variety of equations (hyperbolic and/or parabolic equations) as well as applications such as artificial intelligence. The error analysis above applies and, importantly, we can guarantee the v	alidity of an a priori and {\bf quantitative error bound}. 
In many cases of interest, depending upon the choice of the kernels, 
we can check numerically or theoretically, that the error rate is the optimal convergence rate.

We explored some industrial applications in mathematical finance and non-linear hyperbolic-parabolic equations. The overall algorithm we have developped has been found to be
 robust, fast, accurate and was quite efficient in order to compute standard risk measures for mathematical finance. Indeed, since we can argue that these methods exhibit a sharp convergence rate, they tend to minimize the algorithmic work 	and computational time.

\end{multicols}

\end{document}